\documentclass[10pt,twoside]{amsart}
\usepackage{amssymb,amsbsy,amsmath,amsfonts,amssymb,amscd}
\usepackage{latexsym,euscript,exscale,epic,eepic,epsfig}

\usepackage [frenchb]{babel}

\usepackage{times}
\usepackage{setspace}
\input utile.def

\newtheorem*{theoreme}{Th\'eor\`eme principal}
\newtheorem{proposition}{Proposition}
\newtheorem*{definition}{D\'efinition\rm}
\newtheorem*{remarque}{Remarque}
\newtheorem{exemple}{\it Exemple\/}

\newtheorem{conjecture}{Conjecture\rm}

\begin{document}
\date{\today}
\author{Charles Favre}
\address{Universit\'e Paris 7\\
         UFR de Math\'ematiques\\
         Equipe G\'eom\'etrie et Dynamique\\
         F-75251 Paris Cedex 05\\
         France}
\email{favre@math.jussieu.fr}

\title{%
Les applications monomiales en deux dimensions.
}
%


\begin{abstract}{%
En utilisant le formalisme des vari\'et\'es toriques, on d\'ecrit
comment rendre une application monomiale alg\'ebriquement stable.

\hfill

\noindent {\small ABSTRACT.} 
Using the formalism of toric varieties, we describe how to make a
monomial application algebraically stable.  
}\end{abstract}




\maketitle
\thispagestyle{empty}


\section*{Introduction} Soit $\f : \PP^2 \self$ une application rationnelle 
de l'espace projectif, de degr\'e topologique $ e \ge 1$, et dont
l'action sur le groupe de Picard est donn\'ee par un entier
$\mbox{deg} (\f)$ que l'on supposera toujours plus grand que $2$. On
v\'erifie que la suite des degr\'es $\mbox{deg} (\f^n)$ est
sous-multiplicative ; on notera $\l_1 = \lim_n \mbox{deg
}(\f^n)^{1/n}$ le premier degr\'e dynamique de $\f$ (voir
\cite{russhi}). Les principales propri\'et\'es ergodiques de $\f$
d\'ependent \'etroitement du rapport $\l_1 / e$ (voir \cite{guedj}).
Lorsque $\l_1 > e$ (resp.  $\l_1 < e$), il est ainsi conjectur\'e que
la plupart des points p\'eriodiques sont de type hyperboliques selles
(resp. r\'epulsifs). Il appara\^\i t donc important de comprendre
pr\'ecis\'ement la nature de la suite $\mbox{deg }(\f^n)$. De mani\`ere
g\'en\'erique, cette suite est multiplicative $\mbox{deg }(\f^n)=
\mbox{deg }(\f)^n$, et par suite $\l_1 = \mbox{deg }(\f)$. En
g\'en\'eral cependant, la situation est perturb\'ee par l'existence de
courbes contract\'ees par $\f$ ou l'un de ses it\'er\'es sur un point
d'ind\'etermination. On est donc amen\'e \`a introduire la
d\'efinition suivante (\cite{forsi}).

\begin{definition}
Une application $\f : X \self$ d'une surface rationnelle est dite
Alg\'ebriquement Stable (AS en abr\'eg\'e) si $\cup _{ k \geq 0}
\f^{-k} I(\f)$ est une union d\'enombrable de points, o\`u $I(\f)$
d\'enote l'ensemble d'ind\'etermination de $\f$.
\end{definition}

Il est \'equivalent de supposer que l'action naturelle $\f^*$ de $\f$
sur $\mbox{Pic }(X)$ est composable au sens o\`u $\f^{n*} = \f^{*n}$
pour tout $n \ge 1$. Lorsque $\f$ est AS, $\l_1$ s'identifie au rayon
spectral de $\f^*$ (voir par exemple \cite{dilfa}), le type de
croissance des degr\'es se ram\`ene alors \`a la description du
spectre de $\f^*$. Pour les applications birationnelles de surfaces,
on peut montrer l'existence d'un mod\`ele birationnel de $\PP^2$ dans
lequel $\f$ devient AS, et ainsi d\'ecrire pr\'ecis\'ement la nature
de $\mbox{deg }(\f^n)$ (\cite{dilfa}).  Nous nous proposons de
d\'ecrire des exemples simples d'applications rationnelles pour
lesquels il est impossible de trouver un mod\`ele birationnel dans
lequel celles-ci deviennent AS. Cependant la croissance des degr\'es
de ces applications est compl\`etement \'el\'ementaire \`a d\'ecrire,
et confirme les conjectures suivantes.

\begin{conjecture}\label{conj1}
Soit $\f : \PP^2 \self$ une application rationnelle telle que $ e^2 <
\l_1$. Alors il existe un mod\`ele birationnel $\overline{X}$ de
$\PP^2$, et un rev\^etement ramifi\'e $ h : X \to \overline{X}$ tels
que $\f$ se rel\`eve \`a $X$ et devient AS.
\end{conjecture}

\begin{conjecture}\label{conj2}
Si $\f : \PP^2 \self $ est induit par une application
\emph{polynomiale} de $\CC^2$, il existe une suite d'\'eclatements $\pi
: X \to \PP^2$, telle que $\f$ se rel\`eve \`a $X$ en une application AS.
\end{conjecture}

\begin{conjecture}\label{conj3}
Pour toute application rationnelle $\f : \PP^2 \self $, le premier
degr\'e dynamique est un entier alg\'ebrique.
\end{conjecture}

Soit $ A \= \left[ \begin{array}{cc} a & b \\ c & d \end{array}
\right] \in M ( 2, \ZZ)~, $ une matrice \`a coefficients entiers. On
lui associe l'application \emph{monomiale} $\f _A $ d\'efinie dans
$\CC^2$ par
$$
\f_A (z,w) = (z^a w^b, z^c w^d) ~.
$$
Le degr\'e topologique de $\f_A$ est $| \det A|$, son premier degr\'e
dynamique est le rayon spectral de $A$. Si $|\rho_1| \ge |\rho_2|$
sont les deux valeurs propres de $A$, notons que l'in\'egalit\'e
$\l_1 > e^2$ \'equivaut \`a $|\rho_1| > |\rho_2|$. Lorsque $|\rho_1|
= |\rho_2|$, soit ces deux nombres complexes sont conjugu\'es et non
r\'eels, soit $\rho_1 = \pm \rho_2$.

Nous nous proposons de montrer le
\begin{theoreme}\label{princ}
Soit $\f _A : \PP^2 \self$ une application monomiale de matrice
associ\'ee $A \in M ( 2, \ZZ)$ de d\'eterminant non nul.
\begin{itemize}
\item
Lorsque le spectre de $A$ est r\'eel, on peut trouver une suite
d'\'eclatements $\pi: \overline{X} \to \PP^2$, et un rev\^etement
ramifi\'e $h : X \to \overline{X}$ tels que $\f$ se rel\`eve \`a $X$
et devient AS.
\item
Lorsque le spectre de $A$ est non r\'eel, ses deux valeurs propres
s'\'ecrivent $\exp ( 2 i \pi \theta)$ avec $\theta$ r\'eel (en
particulier $\l_1 = e^2$).

- Si $\theta$ est rationnel, la conclusion pr\'ec\'edente
  s'applique. Plus pr\'ecis\'ement, on peut relever $\f$ en une
  application holomorphe de $\overline{X}$.

- Si $\theta$ est irrationnel, $\f_A$ ne se rel\`eve jamais de
  mani\`ere AS dans un mod\`ele birationnel de $\PP^2$. Plus
  g\'en\'eralement, si on se donne une application surjective $\pi : X
  \to \PP^2$ telle que $\phi$ se rel\`eve \`a $X$ en une application
  rationnelle $\psi$, alors cette derni\`ere n'est pas AS.
\end{itemize}
\end{theoreme}
\begin{remarque}
Lorsque le spectre de $A$ est r\'eel, la preuve du th\'eor\`eme montre
que l'on peut toujours trouver un \emph{morphisme} birationnel $\pi:
X \to \PP^1\times\PP^1$ de telle sorte que
$\pi\circ\f_A\circ\pi^{-1}$ devienne AS.
\end{remarque}

La preuve du th\'eor\`eme s'appuie sur des notions \'el\'ementaires de
g\'eom\'etrie torique (voir \cite{oda}, \cite{fulton}). Nous faisons
tout d'abord quelques rappels succints avant de pr\'esenter la
d\'emonstration de notre r\'esultat. Nous concluons cet article par
quelques exemples.

\section{Surfaces toriques et morphismes toriques}
On se donne un r\'eseau $N \subset \RR^2$ et on note son dual $M \=
\Hom (N, \ZZ)$. On pose $N_{\RR} \= N \otimes \RR$, et $M_{\RR} \= M
\otimes \RR$.  Un \'eventail $\Delta$ est la donn\'ee du singleton $\{
0\}$ et d'une collection finie de semi-droites $\RR_+ v_1, \cdots,
\RR_+ v_s$ avec $v_i \in N$ et de c\^ones $\cC _1 , \cdots, \cC _r$
stricts ($\cC_ i \cap - \cC_i = \{ 0 \}$) de dimension $2$ bord\'es
par des vecteurs dans $\{ v_1, \cdots , v_s\}$ .  On notera 
l'\'eventail ainsi d\'efini
$$
\Delta \= \left\{ \{ 0 \}, \RR_+ v_1, \cdots, \RR_+ v_s , \cC _1 ,
\cdots ,\cC _s \right\} .$$
On associe naturellement \`a tout \'eventail
une surface complexe que l'on note $ \text{T}_N \text{ emb} (\Delta
)$. Une surface de ce type sera dite torique.

Chaque c\^one $\s \in \Delta$ d\'efinit une vari\'et\'e torique affine
comme suit. Consid\'erons le c\^one dual $\check{\s} \= \{ m \in
M_{\RR} , ~ < m , n > \geq 0 \text{ pour tout } n \in \s \}$ et notons
$\cS _{\s} \= \check{\s} \cap M = \{ m \in M , ~ < m , n > \geq 0
\text{ pour tout } n \in \s \}$, l'ensemble des points de ce c\^one
appartenant au r\'eseau $M$.  C'est un semi-groupe pour la loi
additive engendr\'e par un nombre fini d'\'el\'ements par le Lemme de
Gordon. On pose $U_{\s} \= \text{ Spec } \CC [\cS_{\s}]$.  
De mani\`ere concr\`ete, un point $u \in U_{\s}$  est d\'efini par un
morphisme 
$$ u : \cS_{\s} \longrightarrow \CC \text{ tel que } u(0) = 1
\text{ et } u( m + m') = u(m) ~ u (m') \text{ pour tout } m, m' \in
\cS_{\s} .
$$
Si $(m_1, \cdots , m_p)$ est une partie g\'en\'eratrice de $\cS_{\s}$,
l'application de $U_\s$ \`a valeurs dans $\CC^p$ donn\'ee par $u
\mapsto ( u (m_1), \cdots , u (m_p))$ induit un plongement de $U_{\s}$
dans $ \CC ^p$.  Lorsque $\dim ( \s ) = 0, 1$, la surface $U_{\s}$ est
lisse. Lorsque $\dim ( \s ) = 2$, la surface $U_{\s}$ est lisse ssi $
\s = \RR_+ v_1 + \RR_+ v_2$ o\`u $ ( v_1, v_2)$ engendre le r\'eseau
$N$. Lorsque $U_{\s}$ est lisse, on dit que $ \s$ est un c\^one
r\'egulier.

~ \par

Lorsque $\s _1 \subset \s_2$ sont deux c\^ones de $\Delta$, le point $
u \in U_{\s_1}$ induit par restriction sur $\cS _{\s_2} \subset\cS
_{\s_1}$ un point $ u \in U _{\s_2}$. Le morphisme ainsi construit est
un plongement $ U_{\s_1} \hookrightarrow U_{\s_2}$ et $ U_{\s_1} $
forme un ouvert de Zariski dense dans $U_{\s_2}$. En particulier, tout
c\^one contient $\{ 0 \}$ donc $U _{\s}$ contient le tore alg\'ebrique
complexe $ U_{\{ 0 \} }= \Hom ( M , \CC^* ) \cong (\CC ^*)^2$ comme
ouvert dense.  On construit alors la vari\'et\'e $ \text{T}_N \text{
emb} (\Delta ) \= \cup _{\s \in \Delta} U_{\s}$ en recollant chaque
couple $(U_{\s_1} , U_{\s_2})$ le long de $ U_{\s_1 \cap \s_2}$. On
v\'erifie que $ \text{T}_N \text{ emb} (\Delta )$ est compacte ssi le
support de l'\'eventail est $ \RR^2$. On dira alors que $\Delta$ est
complet.  Notons par ailleurs que $ \text{T}_N \text{ emb} (\Delta )$
contient $(\CC^*)^2$ comme ouvert dense et est donc toujours une
surface rationnelle.

~ \par

Les surfaces toriques supportent une action naturelle du tore
alg\'ebrique complexe $$ \text{T}_N \= \Hom ( M , \CC^* ) = N
\otimes_{\ZZ} \CC^* \cong (\CC ^*)^2.$$ Si $t \in \text{T}_N$ et $ u
\in U_{\s} \subset \text{T}_N \text{ emb} (\Delta )$, on pose $(t. u )
(m) \= t (m ) \times u ( m ) \in U_{\s }$.  Notons que le tore
$\text{T}_N $ agit sur $ U_{\{ 0 \} } = \text{T}_N$, par translation.

Il est important de comprendre la d\'ecomposition en orbite de $
\text{T}_N \text{ emb} (\Delta )$ sous l'action de $
\text{T}_N$. Chaque c\^one $\s \in \Delta$ correspond de mani\`ere
biunivoque \`a une orbite $\mbox{orb } (\s)$ de dimension $\dim (\s) +
\dim (\mbox{orb } (\s) ) = 2 $.  Au singleton $\{ 0 \}$, correspond le
tore $ \text{T}_N$; \`a une semi-droite $\RR_+ v$ correspond une
courbe isomorphe \`a $\CC^*$; \`a un c\^one $\cC$ de dimension $2$
correspond un point.

Plus pr\'ecis\'ement, on d\'efinit $ \mbox{orb } (\s) \= \{ u : M \cap
\s^{\perp} \longrightarrow \CC^* \text{ homo. de groupes } \}$. C'est
un tore alg\'ebrique complexe de dimension $ 2 - \dim (\s) $.  Un
point $ u \in \mbox{orb } (\s)$ s'identifie \`a un point de $U_{\s}$
naturellement $\widetilde{u} : \cS_{\s} \longrightarrow \CC$ en posant
$\widetilde{u}(m)= 0$ si $ m \in \cS_{\s} \setminus \s^{\perp}$ et
$\widetilde{u}(m)= u (m) $ si $m \in \s ^{\perp}$. On v\'erifie que
sous cette identification $\mbox{orb } (\s)$ est l'orbite d'un point
de $U_{\s}$ et que les $ \text{T}_N$-orbites sont toutes de cette
forme.

~ \par

D\'ecrivons maintenant les morphismes de surfaces toriques
pr\'eservant cette structure.  Soient $\Delta_1,\Delta_2$ deux
\'eventails d\'efinis respectivement sur des r\'eseaux $N_1, N_2$.
Une application holomorphe $\f: \text{T}_{N_1} \text{ emb} (\Delta_1 )
\to\text{T}_{N_2} \text{ emb} (\Delta_2 )$ est dite torique si elle
pr\'eserve l'action du tore i.e. ssi pour tout $ t \in \text{T}_{N_1},
u \in \text{T}_{N_1} \text{ emb} (\Delta_1 )$, on a
\begin{equation}\label{e1}
\f ( t . u ) = \phi ( t ) . \f (u) ~,
\end{equation}
pour un morphisme de groupe $\phi \in \Hom ( \text{T}_{N_1},
\text{T}_{N_2})$.  Ce morphisme est d\'etermin\'e par une application
$\ZZ$-lin\'eaire $N_1 \to N_2$ que l'on note encore $\phi$.
L'application $\f$ est alors holomorphe ssi pour tout \'el\'ement de
l'\'eventail $\s_1 \in \Delta_1$, on a $\phi ( \s_1 ) \subset \s_2$
pour un $\s_2 \in \Delta_2$.

L'int\'erieur relatif d'un c\^one $\s_i \in \Delta_i$ est par
d\'efinition $\text{Int }( \s_i )= \s_i \setminus \cup \{ \tau \text{
face de } \s \}$.  L'int\'erieur relatif d'un c\^one $\cC$ de
dimension $2$ co\"\i ncide avec son int\'erieur; $\text{Int }( \RR_+ v
) = \RR_+^* v$, et $\text{Int }(\{ 0 \} ) = (\{ 0 \}$.

On v\'erifie alors que $ \f ( \mbox{orb } ( \s) ) = \bigcup \mbox{orb }
( \tau ) $ o\`u $\tau$ parcourt l'ensemble des c\^ones tels que
$\mbox{Int }( \tau ) \cap \phi ( \s ) \not= \emptyset$.  De m\^eme, on
a $ \f^{-1} ( \mbox{orb } ( \s) ) = \bigcup \mbox{orb } ( \tau )$ pour
$\text{Int }( \s ) \cap \phi ( \tau ) \not= \emptyset$.

Notons que le degr\'e topologique de $\f$ est donn\'e par son degr\'e
sur le tore i.e.  $ e = | \det ( \phi ) |$.

Une application m\'eromorphe $\f: \text{T}_{N_1} \text{ emb} (\Delta_1 )
\dashrightarrow \text{T}_{N_2} \text{ emb} (\Delta_2 ) $ 
est dite torique si $ \Gamma$ le graphe de $\f$ est une vari\'et\'e
torique et si les projections naturelles $ \pi_i : \Gamma
\longrightarrow  \text{T}_{N_i} \text{ emb} (\Delta_i )$  pour $i=1,2$
sont des applications toriques. 
Notons $N_{\Gamma}$ le r\'eseau sur lequel l'\'eventail
$\Delta_{\Gamma}$ d\'efinissant $\Gamma$ est construit; et $\phi_i:
N_{\Gamma} \longrightarrow N_i$ pour $i=1,2$ les applications
$\ZZ$-lin\'eaire induisant $\pi_1, \pi_2$.
L'application  $\pi_1$  est une modification propre donc $|\det
(\phi_1)|=1$ et on peut poser $\phi \= \phi_2 \circ \phi_1 ^{-1} : N_1
\to N_2$.
En particulier dans le cas o\`u $\Delta \= \Delta_1= \Delta_2$ et $N
\= N_1 = N_2$ on obtient le fait suivant.

~\par

\emph{ Un morphisme m\'eromorphe $f: \text{T}_{N} \text{ emb}
  (\Delta ) \self$ est d\'etermin\'e de
mani\`ere unique par une application 
$\ZZ$-lin\'eaire de $\phi: N \to N$ i.e. dans une base de $N$ par une
matrice $2 \times 2$ \`a coefficients entiers relatifs.}

~ \par

Notons que l'\'equation \eqref{e1} appliqu\'e au point $u = (1,1)$
montre qu'un morphisme torique $\f$ induit une application monomiale
sur le tore $T_N$.

On d\'eduit ais\'ement des remarques pr\'ec\'edentes:
\begin{itemize}
\item
une courbe $C = \mbox{orb } ( \RR_+ v)$ est contract\'ee sur un point
ssi $\phi (v)$ est dans l'int\'erieur relatif d'un c\^one de $\Delta$
de dimension $2$ ;
\item
un point $p = \mbox{orb } ( \cC)$ est d'ind\'etermination ssi $\phi
(\mbox{Int }\cC)$ contient une semi-droite $\RR_+ v \in \Delta$.
\end{itemize}

Supposons que $\Delta$ soit r\'egulier (i.e. $\text{T}_N \text{ emb}
(\Delta )$ est lisse) et prenons un c\^one de dimension $2$ dans
$\Delta$ donn\'e par $ \cC = \RR_+ v_1 + \RR_+ v_2$.  Soit $v \= v_1 +
v_2$ et consid\'erons $\Delta'$ le raffinement de $\Delta$ obtenu en
rempla\c cant $ \cC$ par les deux c\^ones $\RR_+ v_1 + \RR_+ v$,
$\RR_+ v + \RR_+ v_2$.  L'application lin\'eaire $\id: N \to N$ induit
une modification propre $\pi: \text{T}_N \text{ emb} (\Delta' ) \to
\text{T}_N \text{ emb} (\Delta )$ qui s'identifie \`a l'\'eclatement
de $\text{T}_N \text{ emb} (\Delta )$ au point $\mbox{orb } (\cC )$.

Nous utiliserons un autre exemple de morphismes
toriques. Pour un sous-r\'eseau $N_1 \subset N$ d'indice
fini, l'injection canonique $N_1 \hookrightarrow N$ induit un
rev\^etement $$\phi : \text{T}_{N_1} \text{ emb} (\Delta ) \to
\text{T}_{N} \text{ emb} (\Delta )$$ de degr\'e l'indice $[N_1:N]$.

~ \par

Dans toute la suite, on se place dans le cadre suivant: $\Delta$ est
un \'eventail r\'egulier complet d\'efini sur un r\'eseau $N \subset
\RR ^2$, et $\phi : N \to N$ une application $\ZZ$-lin\'eaire. On
notera toujours $\f$ l'application m\'eromorphe induite
sur $\text{T}_{N} \text{ emb} (\Delta )$. On notera de plus $| \rho
_1| \geq | \rho _2| $ les deux valeurs propres de $\phi$.  On a $\l_1
= | \rho _1|$ et $e = | \rho _1 \times \rho _2| $.  En particulier, on
a toujours $\l_1^2 \geq e$ et $\l_1^2 = e$ ssi $| \rho _1| = | \rho
_2| $.

\section{Le cas $\l_1^2 > e$}

Dans cette section, on suppose que $| \rho _1| > | \rho _2| $ (en
particulier $\rho_1, \rho_2 \in \RR$). Cette in\'egalit\'e est
\'equivalente \`a $\l_1^2 > e$.  On notera $w_1$ et $w_2$ les deux
vecteurs propres associ\'es \`a ces valeurs propres. En g\'en\'eral,
ils ne sont pas \`a coordonn\'ees enti\`eres mais sont d\'efinis sur
une extension quadratique de $\QQ$.

Montrons la
\begin{proposition}{~}\label{prop1}

On peut raffiner $\Delta$ en un \'eventail $\Delta'$ de telle sorte
que $\f : \text{T}_N \text{ emb} (\Delta') \self$ soit AS.  De plus,
il existe un sous-r\'eseau $N' \subset N$ dans lequel $\Delta'$ devient
r\'egulier.
\end{proposition}

Le reste de cette partie est consacr\'ee \`a la preuve de ce r\'esultat.
On posera  $\kappa \= \rho_2 / \rho_1 \in (-1, 1)$.

~ \par

On raffine tout d'abord $\Delta$ en lui ajoutant les semi-droites
$\RR_+ (- v)$ d\`es que $\RR_+ v$ appartient \`a $\Delta$.  L'\'eventail
ainsi obtenu n'est pas n\'ecessairement r\'egulier, mais il est
sym\'etrique par rapport \`a l'origine ce qui permet de raisonner sur
les droites de $\RR^2$ de pentes rationnelles au lieu de semi-droites.

Consid\'erons alors les droites $D_1$ et $D_2$ de $\Delta$ (distinctes
de $\RR w_1$) les plus proches de $\RR w_1$ respectivement dans le
quadrant $\RR_+ w_1 + \RR_+ w_2$ et $\RR_+ w_1 - \RR_+ w_2$. Quitte
\`a raffiner $\Delta$ ces droites existent, et on note $\cC$ l'union
des deux c\^ones convexes qu'elles d\'efinissent. Par construction
celui-ci contient $\RR w_1$.

~ \par

On peut raffiner $\Delta$ de telle sorte que le c\^one $\cC$ s'envoit
strictement dans lui-m\^eme $ \phi (\cC ) \subsetneq \cC$.  En
g\'en\'eral, il n'est pas cependant possible de choisir ce nouvel
\'eventail r\'egulier.

Si $ 1>\kappa > 0$, cette propri\'et\'e est automatiquement
v\'erifi\'ee.  Sinon $-1<\kappa <0$ ce qui \'equivaut \`a $\det \phi <
0$.  Dans la base $(w_1, w_2)$ de $N_{ \RR}$, on \'ecrit $D_1 = \RR
(1, \a)$, $D_2 = \RR ( 1, - \b)$ avec $\a, \b > 0$.  La condition
ci-dessus $ \phi (\cC ) \subsetneq \cC$ est alors \'equivalente \`a
l'in\'egalit\'e
\begin{equation}\label{c1}
|\kappa| \a < \b < |\kappa|^{-1} \a ~.
\end{equation}
On raffine $\Delta$ en lui ajoutant deux droites $\RR ( 1, a')$ et
$\RR (1, -b')$ avec $a' < a$, $b' < b$ de telle sorte que le couple
$(a',b')$ satisfasse \`a \eqref{c1}. Pour ce nouvel \'eventail $ \phi
(\cC ) \subsetneq \cC$.

~\par

De m\^eme, consid\'erons les deux droites $D_1'$ et $D_2'$ de $\Delta$
(distinctes de $\RR w_2$) les plus proches de $\RR w_2$ respectivement
dans le quadrant $\RR_+ w_1 + \RR_+ w_2$ et $\RR_+ w_1 - \RR_+ w_2$,
et $\cD$ l'union des deux c\^ones convexes engendr\'es par celles-ci. Le
m\^eme argument montre que quitte \`a raffiner $\Delta$ on a $ \phi
(\cD ) \supsetneq \cD$.

~ \par

Pour conclure, on ajoute \`a $\Delta$ toutes les droites $D
\not\subset \cD$ telles que $\phi^k ( D ) \in \Delta$ pour un entier
$k \ge 0$. Dans la surface torique associ\'ee \`a cet \'eventail, le
morphisme torique $\f$ induit par $\phi$ poss\`ede au plus deux points
d'ind\'etermination $\mbox{orb } ( \cD_1 )$, $\mbox{orb } ( \cD_2 )$
o\`u $\cD_1, \cD_2$ sont les deux c\^ones convexes d'union $\cD$. Les
courbes contract\'ees par un it\'er\'e de $\f$ sont $\mbox{orb } ( D
)$ pour $D$ une semidroite de $\Delta$, et sont donc \'eventuellement
envoy\'ees sur un des deux points $\mbox{orb } ( \cC_1 )$, $\mbox{orb
} ( \cC_2 )$ o\`u $\cC_1 \cup \cC_2 = \cC$. Ceux-ci sont fixes si
$\det \phi >0$ et permut\'es si $\det \phi < 0$. En g\'en\'eral,
l'\'eventail n'est pas r\'egulier, mais on peut toujours trouver un
sous-r\'eseau de $N$ dans lequel $\Delta$ devient r\'egulier.

Ceci conclut la preuve de la proposition \ref{prop1}.

\section{Le cas $\l_1^2 = e$}

\noindent 
Le cas $|\rho_1| = |\rho_2|$ se s\'epare en 5 sous-cas.

~ \par

$(1)$: $\phi = \rho \id$ avec $\rho = \rho_1 = \rho_2$.  Pour
n'importe quel \'eventail $\Delta$, l'application $ \phi :
\text{T}_{N} \text{ emb} (\Delta) \self$ est holomorphe.

~ \par

$(2)$: $\rho \= \rho_1 = \rho_2$ mais $\phi$ n'est pas diagonalisable.

La droite $D \= \ker ( \phi - \rho \id)$ est rationnelle. Dans ce cas,
on applique l'algorithme pr\'ec\'edent en modifiant la deuxi\`eme
\'etape en rempla\c cant $\cC$ par le c\^one d\'efini par $D$ et le
droite la plus proche de $D$.  On d\'emontre de mani\`ere analogue \`a
la preuve de la proposition \ref{prop1} que le couple $ (\phi, \text{T}_{N}
\text{ emb} (\Delta))$ est AS.

~ \par

$(3)$: $\rho \= \rho_1 = - \rho_2 \ge 1$.

Notons $v_+, v_-$ les vecteurs propres associ\'es resp. \`a $\rho$ et
$ - \rho$. Soit $\sigma$ la reflection fixant $\RR v_+$ et d'axe $\RR
v_-$. Ajoutons \`a $\Delta$ tous les vecteurs $ \pm \RR _+ v, \pm
\RR_+ \sigma (v)$ d\`es que $ \RR _+ v \in \Delta$. Dans ce nouvel
\'eventail $\Delta'$ (singulier en g\'en\'eral), $\phi$ est
holomorphe.  Dans le sous-r\'eseau $N' \= \ZZ v_+ + \ZZ v_-$,
l'\'eventail $\Delta'$ devient r\'egulier.

~ \par

$(4)$: $ \rho_1 = \rho \exp ( i \pi \theta)$, $\rho_2 = \overline
{\rho_1}$ avec $\rho \ge 1$ et $\theta \in \QQ$.

Dans ce cas, $\phi$ est d'ordre fini sur l'ensemble des semi-droites
de $\RR^2$. En ajoutant \`a $\Delta$ toutes les images 
$\phi ^k ( \RR_+ v)$ pour $\RR_+ v \in \Delta$ et $k \geq 0$, on obtient un
raffinement de $\Delta$ dans lequel $\phi$ est holomorphe. 
En g\'en\'eral, l'\'eventail n'est pas r\'egulier.

~ \par

$(5)$: $ \rho_1 = \rho \exp ( i \pi \theta)$, $\rho_2 = \overline
{\rho_1}$ avec $\rho \ge 1$ et $\theta \not\in \QQ$ (dans ce cas on a
en fait toujours $\rho >1$ voir ci-dessous).

Supposons que l'on puisse raffiner $\Delta$ en un \'eventail $\Delta'$
de telle sorte que le couple $ (\phi, \text{T}_{N} \text{ emb} (\Delta'))$ soit
AS. Montrons que ceci ne peut \^etre possible.

Soit $\RR_+ v \in \Delta'$. Comme $\{ \phi ^k (\RR_+ v), k \geq 0 \}$
est dense dans l'ensemble des semi-droites de $\RR^2$, on peut trouver
un entier $k \geq 0$ tel que $\phi ^k (\RR_+ v) \subset \text{Int
} \cC$ pour un c\^one $\cC$ de dimension $2$. Donc toute courbe
torique $\mbox{orb } ( v)$ se contracte sur un point.  De m\^eme si
$\cC \in \Delta'$ est un c\^one de dimension $2$, pour un entier $k
\geq 0$ on a $\text{Int } \cC \supset \RR_+ v \in \Delta'$, et il
s'ensuit que $\mbox{orb } ( \cC) \in \cup _{k \geq 0} I ( \f ^k )$. On
a donc montr\'e que toutes les courbes toriques se contractent sur un
point d'ind\'etermination torique apr\`es it\'eration.

Supposons maintenant que l'on r\'ealise une suite d'\'eclatements
centr\'ee en des points de $\text{T}_{N} \text{ emb} (\Delta)$ qui ne
soient pas toriques (i.e. invariant sous l'action du
tore). L'application induite par $\f$ ne peut lever les
ind\'eterminations toriques, et les courbes toriques contract\'ees
restent critiques car elles sont envoy\'ees sur des points toriques.

On a donc montr\'e que l'on ne peut rendre $ (\phi, \text{T}_{N}
\text{ emb} (\Delta))$ AS en proc\'edant \`a une suite finie
d'\'eclate\-ments toriques ou non.

Supposons maintenant que l'on se donne un morphisme surjectif entre
surfaces lisses $\pi : X \to \text{T}_{N} \text{ emb} (\Delta)$ tel
que $\f$ se rel\`eve \`a $X$ en une application $\psi$.  Notons $U =
\pi^{-1}T_N$, $\cC(\pi)\subset U$ l'ensemble critique de $\pi$ dans
$U$, et $\cE(\pi)\subset \cC(\pi)$ l'ensemble des points o\`u $\pi$
n'est pas localement fini. L'ensemble des valeurs critiques $V_\pi
=\pi\cC(\pi)$ est un sous-ensemble alg\'ebrique propre de $T_N$. Il
n'est pas n\'ecessairement de dimension pure, mais est union
d'une courbe $W_\pi$ et d'un ensemble fini de points $\pi\cE(\pi)=F$.

Soit $x\in V_\pi$, et $z\in\f^{-1}\{x\}$ une pr\'eimage de $x$ par
$\f$. Choisissons $w\in\pi^{-1}\{z\}$, $y\in\pi^{-1}\{x\}$ dans les
fibres respectives au-dessus de $z$ et $y$. Comme $T_N$ est totalement
invariant par $\f$, $w\in \pi^{-1}T_N=U$.

Si $\psi$ est holomorphe en $w$, de l'\'equation
$\pi\circ\psi=\f\circ\pi$ et du fait que $\f$ induit un rev\^etement
non ramifi\'e de $T_N$ sur lui-m\^eme, on d\'eduit que
$z\in\pi\cC(\pi)=V_\pi$. Lorsque $w$ est un point d'ind\'etermination
de $\psi$, son image totale par $\psi$ dans $X$ est une courbe connexe
qui intersecte $U$, et se contracte par $\pi$ sur
$x=\f\circ\pi(w)$. Par suite, $x$ appartient \`a l'ensemble fini de
points $F\subset U$ d\'efini ci-dessus. On a donc d\'emontr\'e que
$$
\f^{-1}V_\pi\subset V_\pi\cup\f^{-1}F~.
$$
Il est alors apparent que la courbe alg\'ebrique $W_\pi\subset V_\pi$
est totalement invariante par $\f$.  On v\'erifie que $\f$ restreint
\`a $T_N$ n'admet aucune courbe totalement invariante. Donc $\pi:
U\setminus\pi^{-1}F\to T_N\setminus F$ est un rev\^etement non
ramifi\'e. En particulier, $U\setminus\pi^{-1}F$ est isomorphe \`a
$(\CC^*)^2\setminus G$ o\`u $G$ est un nombre fini de points, et $\pi$
d\'efinit un morphisme monomial de $(\CC^*)^2\setminus G$ dans
$T_N$. Par extension, $U$ poss\`ede un morphisme birationnel (propre)
sur $(\CC^*)^2$, et est donc obtenue \`a partir du tore par
composition finie d'\'eclatements de points.

La surface $X'$ obtenue par contraction de toute les courbes compactes
de $X$ incluses dans $U$ est une surface alg\'ebrique lisse contenant
le tore $(\CC^*)^2$ comme ouvert dense.  Le morphisme induit par
$\psi$ sur $X'$ pr\'eserve le tore et comme $\pi$ est monomial, $\psi$ agit
aussi comme une application monomiale.

Par contraction successive des courbes exceptionnelles de premi\`ere
esp\`ece de $X'$, on obtient une surface rationnelle minimale $S$
c'est-\`a-dire $\PP^2$, $\PP^1\times\PP^1$ ou une surface de
Hirzebruch. Celles-ci sont toutes des compactifications toriques de
$(\CC^*)^2$.  Les courbes exceptionnelles de premi\`ere esp\`ece de
$X'$ sont toutes hors du tore donc $\psi$ induit aussi un morphisme
monomial sur le tore $(\CC^*)^2\subset S$.

On est donc ramen\'e \`a la situation pr\'ec\'edente o\`u $\pi:X\to S$
est un morphisme birationnel sur une vari\'et\'e torique et $\psi$
induit un morphisme torique dans $S$.  On conclut que $\psi$ ne peut
\^etre AS.  La derni\`ere affirmation du th\'eor\`eme principal est
ainsi d\'emontr\'ee.

\section{Quelques exemples}

Afin de donner un exemple explicite d'application monomiale qui ne
peut \^etre rendu AS dans aucune surface rationnelle, il faut trouver
une matrice $2 \times 2$ \`a coefficients entiers dont les valeurs
propres sont non-r\'eelles, complexes conjugu\'ees $\rho_1, \rho_2 =
\rho \exp ( \pm i \pi \theta)$ avec $\rho >0$ et $\theta \not\in \RR
\setminus \QQ$. Si $A = \left[ \begin{array}{cc} a & b \\ c & d
\end{array} \right]$, le fait d'avoir deux valeurs propres non
r\'eelles est \'equivalent \`a
\begin{equation}\label{e2}
(a-d)^2 + 2 bc < 0 ~.
\end{equation}
En particulier si $a,b,c,d \ge 0$, c'est-\`a-dire lorsque $\f_A$ est
polynomiale dans $\CC^2$, cette in\'egalit\'e n'est jamais
v\'erifi\'ee, et on peut toujours  rendre $\f_A$ AS dans une surface
torique convenable (voir conjecture \ref{conj2}).

De m\^eme lorsque $\f_A$ est birationnel, on a $ad - bc = \pm 1$, et
\eqref{e2} est alors \'equivalent \`a $|\tr A | \le 1$, $\det A = +
1$.  On v\'erifie que dans tous les cas, l'action de $A$ sur les
droites de $\RR^2$ est d'ordre fini, et $\f_A$ est donc holomorphe
dans une surface torique convenable (voir conjecture \ref{conj1}).

~\par

Supposons donc que \eqref{e2} soit v\'erifi\'ee. Si l'argument
$\theta$ des valeurs propres de $A$ est rationnel, on a
\begin{equation} \label{e3}
\theta \in \{0, 1/6, 1/3, 1/4 \} \mod ( 1/2) ~.
\end{equation}
En effet, on remarque que $\rho_1^2 + \rho_2^2 = 2 \rho^2 \cos ( 2 \pi
\theta) = \tr ( \phi^2)$ est un entier, ainsi que $\rho^2 = \det (
\phi)$. Donc $\cos ( 2 \pi \theta), \theta \in \QQ$ ce qui implique
\eqref{e3}.  Les valeurs propres de $A$ sont dans ce cas toujours dans
l'un des corps quadratiques suivants: $\QQ$, $\QQ + \QQ \sqrt{2}$ ou
$\QQ + \QQ \sqrt{3}$.

\begin{exemple}
Pour  $A = \left[ \begin{array}{cc} 1 & 4 \\ -1 & 0
\end{array} \right]$, on a $\rho_1, \rho_2 = 1/2 ( 1 \pm i \sqrt{7})$.
L'application $\f_A (x,y) = (x y^{-1}, x^4)$ ne peut \^etre rendu AS
dans un mod\`ele birationnel de $\PP^2$.
\end{exemple}
Pour conclure mentionnons l'exemple suivant, qui montre qu'en
g\'en\'eral on ne peut rendre une application AS simplement en
\'eclatant, mais qu'il est n\'ecessaire de passer \`a un rev\^etement
ramifi\'e.
\begin{exemple}
Prenons comme \'eventail de d\'epart $$\Delta = \{ \{ 0\}, \pm \RR_+
(0,1), \pm \RR_+ (1,0), \pm \RR_+ (0,1) \pm \RR_+ (1,0) \}~,$$
i.e. $\text{T}_{N} \text{ emb} (\Delta) = \PP^1 \times \PP^1$.  Soit
$$
A = \left[ \begin{array}{cc} 0 & -8 \\ 1 & 4
\end{array} \right] ~.$$
Alors, on ne peut raffiner $\Delta$
en un \'eventail $\Delta'$ \emph{r\'egulier} tel que $\phi$ devienne
AS.

En effet supposons que ce soit le cas. Comme $A$ est d'ordre fini
d'ordre $8$ sur les semi-droites de $\RR^2$, $\f_A$ est
n\'ecessairement holomorphe dans $\text{T}_{N} \text{ emb} (\Delta')$.
Notons $v = (-1, k)$ avec $ k \in \ZZ_+$ le vecteur le plus proche de
$(0,1)$ dans le c\^one $\RR_+ (0,1) + \RR_+ (-1,0)$.  Pour un vecteur
$w \in \ZZ^2$, on note $\widetilde{w}$ le vecteur dans $\RR_+w \cap
\ZZ^2$ de plus petite norme.  Comme $\Delta'$ est r\'egulier les
vecteurs $\widetilde{\phi (0,1)}, \widetilde{\phi (v)}$ forme une base
de $ \ZZ^2$. de mani\`ere \'equivalente, on a
$$|\det (\widetilde{\phi (0,1)}, \widetilde{\phi (v)} )| = 1~.$$ Or
$\widetilde{\phi (0,1)} = ( -2, 1 )$, et $\widetilde{\phi (v)} = (-
8k, 4k-1 )/ (8k \wedge 4k-1 )$. On v\'erifie que ceci aboutit \`a une
contradiction.
\end{exemple}

{\bf Remerciements:} je remercie vivement V.~Guedj pour les discussions que
nous avons eues sur ce travail.


%
\end{document}